\documentclass[10pt,twoside,reqno]{amsart}
\usepackage{amssymb}
\usepackage{multirow}
\calclayout

\numberwithin{equation}{section}
\makeatletter

\renewcommand{\@secnumfont}{\bfseries}

\renewcommand{\section}{\@startsection{section}{1}%
  {0mm}{.7\linespacing\@plus\linespacing}{.5\linespacing}
  {\normalfont\bfseries\centering}}

\newcommand{\bibsection}{\@startsection{section}{1}%
  {0mm}{.7\linespacing\@plus\linespacing}{.5\linespacing}
  {\normalfont\scshape\centering}}

\renewcommand{\@biblabel}[1]{#1.}

\newtheorem{thm}{\bf Theorem}[section]

\begin{document}

\vspace{1.3cm}

\title {Extended degenerate Stirling numbers of the second kind and extended degenerate Bell polynomials}

\author{Taekyun Kim}
\address{Department of Mathematics, College of Science, Tianjin Polytechnic University, Tianjin City, 300387, China.}
\address{Department of Mathematics, Kwangwoon University, Seoul 139-701, Republic
	of Korea}
\email{tkkim@kw.ac.kr}

\author{Dae San Kim}
\address{Department of Mathematics, Sogang University, Seoul 121-742, Republic of Korea}
\email{dskim@sogang.ac.kr}


\begin{abstract} 
In a recent work, the degenerate Stirling polynomials of the second kind were studied by T. Kim. In this paper, we investigate the extended degenerate Stirling numbers of the second kind and the extended degenerate Bell polynomials associated with them. As results, we give some expressions, identities and properties about the extended degenerate Stirling numbers of the second kind and the extended degenerate Bell polynomials.
\end{abstract}
\maketitle

\markboth{\centerline{\scriptsize Extended degenerate Stirling numbers of the second kind and Bell polynomials }}
{\centerline{\scriptsize T. Kim, D. S. Kim }}

\bigskip
\medskip
\section{Introduction}
As is well known, the Stirling numbers of the first kind are defined by the falling factorial sequence as
\begin{equation}\begin{split}\label{01}
(x)_n = \sum_{k=0}^n S_1(n,k) x^k, \,\,(n \geq 0),\quad (\textnormal{see} \,\, [1-16]),
\end{split}\end{equation}
where $(x)_0=1$, $(x)_n = x(x-1)\cdots(x-n+1)$, $(n \geq 1)$. The Stirling numbers of the second kind are defined by the generating function 
\begin{equation}\begin{split}\label{02}
\frac{1}{k!}(e^t-1)^k = \sum_{n=k}^\infty S_2(n,k) \frac{t^n}{n!},\quad (\textnormal{see} \,\, [1-16]).
\end{split}\end{equation}
From \eqref{02}, we note that
\begin{equation}\begin{split}\label{03}
x^n = \sum_{k=0}^n S_2(n,k) (x)_k,\,\,(n \geq 0),\quad (\textnormal{see} \,\, [6-8,11,12]).
\end{split}\end{equation}

In [2,3], L. Carlitz introduced the degenerate Stirling, Bernoulli and Eulerian numbers. With the viewpoint of  generalizing Stirling numbers, the $r$-Stirling numbers of the first kind and of the second kind were introduced by Broder (see [1]). It is known that the $r$-Stirling numbers of the second kind are given by the generating function 

\begin{equation}\begin{split}\label{04}
e^{rt} \frac{1}{k!} (e^t-1)^k = \sum_{n=0}^\infty S_{2,r} (n+r,k+r) \frac{t^n}{n!},\quad (\textnormal{see} \,\, [8]).
\end{split}\end{equation}
The Bell polynomials are defined by 
\begin{equation}\begin{split}\label{05}
Bel_n(x) = \sum_{k=0}^n S_2(n,k) x^k,\quad (\textnormal{see} \,\, [1,8]).
\end{split}\end{equation}
Thus, by \eqref{02} and \eqref{05}, we get

\begin{equation}\begin{split}\label{06}
e^{x(e^t-1)} = \sum_{n=0}^\infty Bel_n(x) \frac{t^n}{n!},\quad (\textnormal{see} \,\, [11]).
\end{split}\end{equation}
When $x=1$, $Bel_n=Bel_n(1)$, $(n \geq 0)$, are called the Bell numbers.

In [2], L. Carlitz introduced the degenerate factorial sequences given by
\begin{equation}\begin{split}\label{07}
(x|\lambda )_0 = 1,\,\, (x|\lambda )_n = x(x-\lambda )(x-2\lambda )\cdots(x-(n-1)\lambda ),\,\, (n \geq 1),
\end{split}\end{equation}
where $\lambda  \in \mathbb{R}$. Note that $\lim_{\lambda  \rightarrow 0}(x|\lambda )_n = x^n$
and $\lim_{\lambda  \rightarrow 1}(x|\lambda )_n = (x)_n$, $(n \geq 0)$.

Recently, the degenerate Stirling numbers of the second kind are defined by the generating function 
\begin{equation}\begin{split}\label{08}
\frac{1}{k!} ((1+\lambda t)^{\frac{1}{\lambda }}-1)^k = \sum_{n=k}^\infty S_{2,\lambda }(n,k) \frac{t^n}{n!},
\end{split}\end{equation}
where $k \in \mathbb{N} \cup \{0\}$ and $\lambda \in \mathbb{R}$ (see [8]). Note that $\lim_{\lambda  \rightarrow 0}S_{2,\lambda }(n,k) = S_2(n,k)$, $(n,k \geq 0)$.

In [11], Kim-Kim-Dolgy introduced the degenerate Bell polynomials associated with the degenerate Stirling numbers of the second kind given by the generating function 
\begin{equation}\begin{split}\label{09}
e^{x((1+\lambda t)^{\frac{1}{\lambda }}-1)} = \sum_{n=0}^\infty Bel_{n,\lambda }(x) \frac{t^n}{n!}.
\end{split}\end{equation}
where $x=1$, $Bel_{n,\lambda }=Bel_{n,\lambda }(1) $ are called the degenerate Bell numbers.

From \eqref{09}, we have

\begin{equation}\begin{split}\label{10}
Bel_{n,\lambda }(x) = \frac{1}{e^x} \sum_{k=0}^\infty \frac{1}{k!} (k|\lambda )_n x^k, \,\,(n \geq 0),
\end{split}\end{equation} 
\begin{equation}\begin{split}\label{11}
Bel_{n,\lambda }(x) = x \sum_{k=1}^n \sum_{j=1}^k {k-1 \choose j-1} S_1(n,k) \lambda ^{n-k} Bel_{j-1}(x),
\end{split}\end{equation}
and
\begin{equation}\begin{split}\label{12}
Bel_{n,\lambda }(x) = \sum_{k=0}^n \sum_{m=0}^k S_2(k,m) S_1(n,k) \lambda ^{n-k} x^m,\quad (\textnormal{see} \,\, [11]).
\end{split}\end{equation}

In this paper, we investigate the extended degenerate Stirling numbers of the second kind and the extended degenerate Bell polynomials associated with them. As results, we give some expressions, identities and properties about the extended degenerate Stirling numbers of the second kind and the extended degenerate Bell polynomials.

\section{Extended degenerate Stirling numbers of the second kind and extended degenerate Bell polynomials}

Throughout this section, we assume that $\lambda \in \mathbb{R}$.

From \eqref{03}, we note that
\begin{equation}\begin{split}\label{13}
&\sum_{n=0}^\infty Bel_{n,\lambda }(x) \frac{t^n}{n!} = \sum_{k=0}^\infty x^k \frac{1}{k!} ((1+\lambda t)^{\frac{1}{\lambda }}-1)^k\\
&= \sum_{k=0}^\infty x^k \sum_{n=k}^\infty S_{2,\lambda }(n,k) \frac{t^n}{n!} = \sum_{n=0}^\infty \left( \sum_{k=0}^n x^k S_{2,\lambda }(n,k) \right) \frac{t^n}{n!}.
\end{split}\end{equation}
By comparing the coefficients on both sides of \eqref{13}, we get
\begin{equation*}\begin{split}
Bel_{n,\lambda }(x) = \sum_{k=0}^n x^k S_{2,\lambda }(n,k), \,\,(n \geq 0).
\end{split}\end{equation*}
For $r \in \mathbb{N} \cup \{0\}$, we define the extended degenerate Stirling numbers of the second kind as
\begin{equation}\begin{split}\label{14}
\frac{1}{k!} (1+\lambda t)^{\frac{r}{\lambda }} ((1+\lambda t)^{\frac{1}{\lambda }}-1)^k = \sum_{n=k}^\infty S_{2,r}(n+r,k+r|\lambda ) \frac{t^n}{n!}.
\end{split}\end{equation}
From \eqref{14}, we note that
\begin{equation}\begin{split}\label{15}
&\frac{1}{k!}((1+\lambda t)^{\frac{1}{\lambda }}-1)^k (1+\lambda t)^{\frac{r}{\lambda }}= \left( \sum_{l=k}^\infty S_{2,\lambda }(l,k) \frac{t^l}{l!} \right) \left( \sum_{m=0}^\infty \left(\frac{r}{\lambda }\right)^m \frac{(\log(1+\lambda t)^m}{m!} \right)\\
&= \left( \sum_{l=k}^\infty S_{2,\lambda }(l,k) \frac{t^l}{l!} \right) \left( \sum_{m=0}^\infty r^m \lambda ^{-m} \sum_{i=m}^\infty S_1 (i,m) \lambda ^i \frac{t^i}{i!} \right)\\
&= \left( \sum_{l=k}^\infty S_{2,\lambda }(l,k) \frac{t^l}{l!} \right) \left( \sum_{i=0}^\infty \left( \sum_{m=0}^i r^m \lambda ^{i-m} S_1(i,m) \right) \frac{t^i}{i!} \right)\\
&= \sum_{n=k}^\infty \left( \sum_{l=k}^n \sum_{m=0}^{n-l} {n \choose l} r^m \lambda ^{n-m-l} S_1(n-l,m) S_{2,\lambda }(l,k) \right) \frac{t^n}{n!}
\end{split}\end{equation}
Therefore, by \eqref{14} and \eqref{15}, we obtain the following theorem.

\begin{thm}
For $n,k \in \mathbb{N} \cup \{0\}$ with $n \geq k$, we have
\begin{equation*}\begin{split}
S_{2,r}(n+r,k+r|\lambda ) = \sum_{l=k}^n \sum_{m=0}^{n-l} {n \choose l} r^m \lambda ^{n-m-l} S_1(n-l,m) S_{2,\lambda }(l,k).
\end{split}\end{equation*}
\end{thm}

Note that
\begin{equation*}\begin{split}
\lim_{\lambda  \rightarrow 0}S_{2,r}(n+r,k+r|\lambda ) &= \sum_{l=k}^n {n \choose l} r^{n-l} \lim_{\lambda  \rightarrow 0} S_{2,\lambda }(l,k)\\
&= \sum_{l=k}^n {n \choose l} r^{n-l} S_2(l,k) = S_{2,r}(n+r,k+r).
\end{split}\end{equation*}
On the other hand,
\begin{equation}\begin{split}\label{16}
& \frac{1}{k!} (1+\lambda t)^{\frac{r}{\lambda }} ((1+\lambda t)^{\frac{1}{\lambda }}-1)^k = \frac{1}{k!} ((1+\lambda t)^{\frac{1}{\lambda }}-1+1)^r (1+\lambda t)^{\frac{1}{\lambda }}-1)^k\\
&= \frac{1}{k!} \sum_{m=0}^\infty {r \choose m} ((1+\lambda t)^{\frac{1}{\lambda }}-1)^{m+k} = \sum_{m=0}^\infty {r \choose m} \frac{	(m+k)!}{	k!} \cdot \frac{1}{(m+k)!} ((1+\lambda t)^{\frac{1}{\lambda }}-1)^{m+k}\\
&=\sum_{m=0}^\infty {r \choose m} m! {m+k \choose m} \sum_{n=m+k}^\infty S_{2,\lambda }(n,m+k) \frac{t^n}{n!} \\
&= \sum_{n=k}^\infty \left( \sum_{m=0}^{n-k} {r \choose m} m! {m+k \choose m} S_{2,\lambda }(n,m+k) \right) \frac{t^n}{n!}
 \end{split}\end{equation}
Thus, by \eqref{14} and \eqref{16}, we get
\begin{equation}\begin{split}\label{17}
S_{2,r}(n+r,k+r|\lambda ) = \sum_{m=0}^{n-k} {m+k \choose m}{r \choose m} m!  S_{2,\lambda }(n,m+k),
\end{split}\end{equation}
where $n,k \geq 0$ with $n \geq k$.

Now, we observe that
\begin{equation}\begin{split}\label{18}
(1+\lambda t)^{\frac{x+r}{\lambda }}& = (1+\lambda t)^{\frac{r}{\lambda }} (1+\lambda t)^{\frac{x}{\lambda }} = (1+\lambda t)^{\frac{r}{\lambda }} ((1+\lambda t)^{\frac{1}{\lambda }}-1+1)^x\\
&= (1+\lambda t)^{\frac{r}{\lambda }}\sum_{k=0}^\infty (x)_k \frac{1}{k!} ((1+\lambda t)^{\frac{1}{\lambda }}-1)^k\\
&= \sum_{k=0}^\infty (x)_k \frac{1}{k!} ((1+\lambda t)^{\frac{1}{\lambda }}-1)^k (1+\lambda t)^{\frac{r}{\lambda }} \\
&= \sum_{k=0}^\infty (x)_k \sum_{n=k}^\infty S_{2,r}(n+r,k+r|\lambda ) \frac{t^n}{n!}\\
&= \sum_{n=0}^\infty \left( \sum_{k=0}^n (x)_k S_{2,r}(n+r,k+r|\lambda ) \right) \frac{t^n}{n!},
\end{split}\end{equation}
and
\begin{equation}\begin{split}\label{19}
(1+\lambda t)^{\frac{x+r}{\lambda }} = \sum_{n=0}^\infty \left( \frac{x+r}{\lambda } \right)_n \lambda ^n \frac{t^n}{n!} = \sum_{n=0}^\infty (x+r|\lambda )_n \frac{t^n}{n!}.
\end{split}\end{equation}
On the other hand,
\begin{equation}\begin{split}\label{20}
(1+\lambda t)^{\frac{x+r}{\lambda }}&= e^{\frac{x+r}{\lambda }\log(1+\lambda t)} = \sum_{k=0}^\infty \left( \frac{x+r}{\lambda} \right)^k \frac{1}{k!} (\log(1+\lambda t))^k\\
&= \sum_{k=0}^\infty (x+r)^k \lambda ^{-k} \sum_{n=k}^\infty S_1(n,k) \lambda ^n \frac{t^n}{n!}\\
&= \sum_{n=0}^\infty \left( \sum_{k=0}^n \lambda ^{n-k} (x+r)^k S_1(n,k) \right) \frac{t^n}{n!}.
\end{split}\end{equation}

Therefore, by \eqref{18},\eqref{19} and \eqref{20}, we obtain the following theorem.

\begin{thm}
For $n \geq 0$, we have
\begin{equation*}\begin{split}
(x+r|\lambda )_n &= \sum_{k=0}^n S_{2,r}(n+r,k+r|\lambda )(x)_k\\
&= \sum_{k=0}^n \lambda ^{n-k} S_1(n,k) (x+r)^k.
\end{split}\end{equation*}
\end{thm}

In particular,
\begin{equation*}\begin{split}
S_{2,r}(n+r,k+r|\lambda ) = \sum_{m=0}^{n-k} {m+k \choose m} {r \choose m} m! S_{2,\lambda }(n,m+k),
\end{split}\end{equation*}
where $n,k \geq 0$ with $n \geq k$.

In view of \eqref{06}, we define the extended degenerate Bell polynomials associated with the extended degenerate Stirling numbers of the second kind as follows:
\begin{equation}\begin{split}\label{21}
(1+\lambda t)^{\frac{r}{\lambda }}e^{x((1+\lambda t)^{\frac{1}{\lambda }}-1)} = \sum_{n=0}^\infty Bel_{n,\lambda }^{(r)} (x) \frac{t^n}{n!}.
\end{split}\end{equation}
From \eqref{21}, we note that
\begin{equation}\begin{split}\label{22}
(1+\lambda t)^{\frac{r}{\lambda }}e^{x((1+\lambda t)^{\frac{1}{\lambda }}-1)}&= \sum_{k=0}^\infty x^k \frac{1}{k!} ((1+\lambda t)^{\frac{1}{\lambda }}-1)^k (1+\lambda t)^{\frac{r}{\lambda }} \\
&= \sum_{k=0}^\infty x^k \sum_{n=k}^\infty S_{2,r} (n+r,k+r|\lambda ) \frac{t^n}{n!}\\
&= \sum_{n=0}^\infty \left( \sum_{k=0}^n x^k S_{2,r}(n+r,k+r|\lambda ) \right) \frac{t^n}{n!}.
\end{split}\end{equation}
Therefore, by \eqref{21} and \eqref{22}, we obtain the following theorem.

\begin{thm}
For $n \geq 0$, we have
\begin{equation*}\begin{split}
Bel_{n,\lambda }^{(r)}(x) &= \sum_{k=0}^n x^k S_{2,r}(n+r,k+r|\lambda )\\
&= \sum_{k=0}^n \left( \sum_{m=k}^n \sum_{l=0}^{n-m} r^l \lambda ^{n-m-l} S_1(n-m,l) S_{2,\lambda }(m,k) {n \choose m} \right) x^k.
\end{split}\end{equation*}	
\end{thm}

Remark. When $x=1$, $Bel_{n,\lambda }^{(r)} = Bel_{n,\lambda }^{(r)}(1)$ are called the extended degenerate Bell numbers associated with the extended degenerate Stirling numbers of the second kind.
\begin{equation*}\begin{split}
Bel_{n,\lambda }^{(r)} &= \sum_{k=0}^n S_{2,r} (n+r,k+r|\lambda )\\
&= \sum_{k=0}^n \sum_{m=k}^n \sum_{l=0}^{n-m} {n \choose m} r^l \lambda ^{n-m-l} S_1(n-m,l) S_{2,\lambda }(m,k).
\end{split}\end{equation*}
From \eqref{14}, we note that

\begin{equation}\begin{split}\label{23}
& \frac{1}{k!} (1+\lambda t)^{\frac{r}{\lambda }}((1+\lambda t)^{\frac{1}{\lambda }}-1)^k = \frac{1}{k!} \sum_{l=0}^k {k \choose l} (-1)^{k-l} (1+\lambda t)^{\frac{l+r}{\lambda }}\\
&= \frac{1}{k!} \sum_{l=0}^k {k \choose l} (-1)^{k-l} e^{\frac{l+r}{\lambda }\log(1+\lambda t)} \\&= \frac{1}{k!} \sum_{l=0}^k {k \choose l} (-1)^{k-l} \sum_{m=0}^\infty \left( \frac{l+r}{\lambda } \right)^m \frac{(\log(1+\lambda t))^m}{m!}\\
&=\frac{1}{k!}\sum_{l=0}^k {k \choose l} (-1)^{k-l} \sum_{m=0}^\infty \lambda ^{-m} (l+r)^m \sum_{n=m}^\infty S_1(n,m) \lambda ^n \frac{t^n}{n!}\\
&= \frac{1}{k!} \sum_{l=0}^\infty {k \choose l} (-1)^{k-l} \sum_{n=0}^\infty \sum_{m=0}^n \lambda ^{n-m} (l+r)^m S_1(n,m) \frac{t^n}{n!}\\
&= \sum_{n=0}^\infty \left\{ \frac{1}{k!} \sum_{m=0}^n \lambda ^{n-m} S_1(n,m) \sum_{l=0}^k {k \choose l} (-1)^{k-l} (l+r)^m \right\} \frac{t^n}{n!}\\
&= \sum_{n=0}^\infty \left\{ \frac{1}{k!} \sum_{m=0}^n \lambda ^{n-m} S_1(n,m) \Delta^k r^m \right\} \frac{t^n}{n!},
\end{split}\end{equation}
where $\Delta f(x) = f(x+1) - f(x)$.

Therefore, by \eqref{14} and \eqref{23}, we obtain the following theorem.

\begin{thm}
For $n,k \geq 0$, we have
\begin{equation*}\begin{split}
\frac{1}{k!} \sum_{m=0}^n \lambda^{n-m} S_1(n,m) \Delta^k r^m = \begin{cases}
0&\text{if}\,\, n <k    \\
S_{2,r}(n+r,k+r|\lambda)&\text{if}\,\,n \geq k.
\end{cases}
\end{split}\end{equation*}
\end{thm}
From Theorem 2.3 and Theorem 2.4, we have

\begin{equation}\begin{split}\label{24}
Bel_{n,\lambda }^{(r)}(x)& = \sum_{k=0}^n x^k \frac{1}{k!} \sum_{m=0}^n \lambda ^{n-m} S_1(n,m) \Delta^k r^m\\
&= \sum_{m=0}^n \lambda ^{n-m} S_1(n,m) \sum_{k=0}^n x^k \frac{1}{k!} \Delta^k r^m.
\end{split}\end{equation}

By \eqref{21}, we get

\begin{equation}\begin{split}\label{25}
&\sum_{n=0}^\infty Bel_{n,\lambda }^{(r)}(x) \frac{t^n}{n!} = (1+\lambda t)^{\frac{r}{\lambda }} e^{x ((1+\lambda t)^{\frac{1}{\lambda }}-1)}\\
&= \left( \sum_{l=0}^\infty Bel_{l,\lambda }(x) \frac{t^l}{l!} \right) \left( \sum_{m=0}^\infty \lambda ^{-m} r^m \frac{1}{m!} (\log(1+\lambda t))^m \right)\\
&= \left( \sum_{l=0}^\infty Bel_{l,\lambda }(x) \frac{t^l}{l!} \right) \left( \sum_{k=0}^\infty \left( \sum_{m=0}^k \lambda ^{k-m} r^m S_1(k,m) \right) \frac{t^k}{k!} \right)\\
&= \sum_{n=0}^\infty \left( \sum_{k=0}^n \sum_{m=0}^k {n \choose k} Bel_{n-k,\lambda } (x) \lambda ^{k-m} r^m S_1(k,m) \right) \frac{t^n}{n!}.
\end{split}\end{equation}

Therefore, by comparing the coefficients on both sides of \eqref{25}, we obtain the following theorem.

\begin{thm}
For $n \geq 0$, we have
\begin{equation*}\begin{split}
Bel_{n,\lambda }^{(r)}(x) = \sum_{k=0}^n \sum_{m=0}^k {n \choose k} Bel_{n-k,\lambda } (x) \lambda ^{k-m} r^m S_1(k,m).
\end{split}\end{equation*}
\end{thm}

From \eqref{21}, we have
\begin{equation}\begin{split}\label{26}
&\sum_{n=0}^\infty Bel_{n,\lambda }^{(r)}(x) \frac{t^n}{n!} = (1+\lambda t)^{\frac{r}{\lambda }}e^{x((1+\lambda t)^{\frac{1}{\lambda }}-1)} \\
&= e^{-x}\sum_{k=0}^\infty \frac{x^k}{k!} (1+\lambda t)^{\frac{k+r}{\lambda }} \\
&= e^{-x} \sum_{k=0}^\infty \frac{x^k}{k!} \sum_{m=0}^\infty \left( \frac{k+r}{\lambda }\right)^m \frac{1}{m!} (\log(1+\lambda t))^m \\
&= e^{-x} \sum_{k=0}^\infty \frac{x^k}{k!} \sum_{n=0}^\infty \left( \sum_{m=0}^n \lambda ^{n-m} (k+r)^m S_1(n,m) \right) \frac{t^n}{n!} \\
&= e^{-x}\sum_{n=0}^\infty \left( \sum_{m=0}^n \lambda ^{n-m} S_1(n,m) \sum_{k=0}^\infty \frac{x^k}{k!} (k+r)^m \right) \frac{t^n}{n!}.
\end{split}\end{equation}
Comparing the coefficients on both sides of \eqref{26}, we have
\begin{equation}\begin{split}\label{27}
Bel_{n,\lambda }^{(r)}(x) =  e^{-x}\sum_{m=0}^n \lambda ^{n-m} S_1(n,m) \sum_{k=0}^\infty \frac{x^k}{k!} (k+r)^m,
\end{split}\end{equation}
where $n \geq 0$.

From \eqref{14}, we have
\begin{equation}\begin{split}\label{28}
&\frac{1}{m!} (1+\lambda t)^{\frac{r}{\lambda }} ((1+\lambda t)^{\frac{1}{\lambda }}-1)^m \frac{1}{k!} ((1+\lambda t)^{\frac{1}{\lambda }}-1)^k \\
&= \frac{1}{m!k!} (1+\lambda t)^{\frac{r}{\lambda }} ((1+\lambda t)^{\frac{1}{\lambda }}-1)^{m+k} = \frac{(m+k)!}{m!k!}(1+\lambda t)^{\frac{r}{\lambda}}\frac{((1+\lambda t)^{\frac{1}{\lambda }}-1)^{m+k} }{(m+k)!}\\
&= {m+k \choose m} \sum_{n=m+k}^\infty S_{2,r}(n+r,m+k+r|\lambda ) \frac{t^n}{n!}.
\end{split}\end{equation}

On the other hand,

\begin{equation}\begin{split}\label{29}
&\frac{1}{m!} (1+\lambda t)^{\frac{r}{\lambda }} ((1+\lambda t)^{\frac{1}{\lambda }}-1)^m \frac{1}{k!} ((1+\lambda t)^{\frac{1}{\lambda }}-1)^k \\
&=\left( \sum_{l=m}^\infty S_{2,r} (l+r,m+r|\lambda ) \frac{t^l}{l!} \right) \left( \sum_{j=k}^\infty S_{2,\lambda}(j,k) \frac{t^j}{j!} \right)\\
&= \sum_{n=k+m}^\infty \left( \sum_{l=m}^{n-k} {n \choose l} S_{2,r}(l+r,m+r|\lambda ) S_{2,\lambda }(n-l,k) \right) \frac{t^n}{n!}
\end{split}\end{equation}

Therefore, by \eqref{28} and \eqref{29}, we obtain the following theorem.

\begin{thm}
For $n,m,k \geq 0$ with $n \geq m+k$, we have
\begin{equation*}\begin{split}
 {m+k \choose m} S_{2,r}(n+r,m+k+r|\lambda ) = \sum_{l=m}^{n-k} {n \choose l} S_{2,r}(l+r,m+r|\lambda ) S_{2,\lambda }(n-l,k).
\end{split}\end{equation*}
\end{thm}

\end{document}